\documentclass[11pt]{article}
\usepackage{amsmath}
\usepackage{amssymb}
\usepackage{amsthm}

\newtheorem{theo}{Theorem}[section]
\newtheorem{lm}{Lemma}[section]
\newtheorem{df}{Definition}[section]
\newtheorem{ex}{Example}[section]
\newtheorem{cor}{Corollary}[section]
\newtheorem{rmk}{Remark}[section]
\newtheorem{cond}{Condition}[section]

\allowdisplaybreaks

\numberwithin{equation}{section}

\def\R{{\mathbb R}}
\def\Z{{\mathbb Z}}
\def\N{{\mathbb N}}
\newcommand\M{{\mathfrak M}}

\def\E{{\mathbf E}}
\def\P{{\mathbf P}}
\def\1{{\mathbf 1}}

\def\const{\mathop{\rm const}}
\def\aff{\mathop{{\rm aff}}}
\def\ri{\mathop{{\rm ri}}}
\def\cl{\mathop{{\rm cl}}}
\def\eps{\varepsilon}
\let\phi=\varphi
\def\qed{\hfill\rule{.2cm}{.2cm}}
\newcommand{\K}{{\kappa_1+\dotsb+\kappa_K}}

\allowdisplaybreaks
\title{Introduction to shape stability for a storage model}
\author{M.V.~Menshikov$^{1}$ 
\and
 V.V.~Sisko$^{2}$
\and M.~Vachkovskaia$^{3}$}

\begin{document}
\maketitle
{\footnotesize
\noindent
$^1$University of Durham,
Department of Mathematical Sciences,
South Road, Durham DH1 3LE, UK.\\
E-mail: Mikhail.Menshikov@durham.ac.uk

\noindent
$^2$ Institute of Mathematics and Statistics,
Federal University of Fluminense,
R. M\'ario Santos Braga, s/n,
CEP 24020-140 Niter\'oi, RJ, Brazil.\\
E-mail: valentin@mat.uff.br

\noindent
$^3$ Department of Statistics, Institute of Mathematics, Statistics and Scientific Computation, University of Campinas -- UNICAMP, rua S\'ergio Buarque de Holanda 651, CEP
13083--859, Campinas SP, Brazil\\
E-mail: marinav@ime.unicamp.br

}

\begin{abstract}
  We consider a new idea for a storage model on $n$ nodes, namely stability of  shape. 
These nodes support   $K $ neighborhoods $S_i\subset \{1,
  \ldots, n\}$ and items arrive at the $S_i$ as independent Poisson
  streams with rates $\lambda_i$, $i=1$, $\ldots\,$, $K $. Upon
  arrival at $S_i$ an item is stored at node $j \in S_i$ where $j$ is
  determined by some policy. 
Under natural conditions on the $\lambda_i$ we
  exhibit simple local policies such that the multidimensional process describing the evolution 
of the number of items at each node is positive recurrent
  (stable) in shape.
\end{abstract}



\section{Description of the model}
Stability in shape is of interest in several models. There are of
course various growth models, see for example the crystal growth model
studied in~\cite{AMS}, though the methods used there are very
different from those we use in this paper. Another model which is
relevant is a queueing system with server vacations or maintenance
periods where stability in shape can be seen as a fairness criterion
for arriving jobs. 
It is also reasonable to view our storage model as a simplified
version of the supermarket model (by dropping the service), see for example~\cite{LMcD}. 

We have chosen to focus on the routing aspect of the model
here. Rather more complex phenomena appear when service is considered
as well and we are investigating a model in which
service times are dependent upon both the arrival neighborhood and
the allocated server.

We consider a storage system (or library) with a finite number of
nodes where identical items are to be stored. The $n$ nodes support
non-empty neighborhoods $S_i$, $i=1$, $\ldots\,$, $K$ with 
 \[\bigcup_{i=1}^K  S_i=\{1, \ldots, n\},\]
 and $1\le K \le 2^n-1$. Items
arrive at the neighborhoods as independent Poisson processes with
rates $\lambda_i > 0$ at $S_i$, $i=1$, $\ldots\,$, $K$ where we
suppose that $\sum_{i=1}^n \lambda_i = 1$. Upon arrival at $S_i$ an
item is stored at a node $j \in S_i$ where $j$ is chosen by some
policy. We consider \emph{local} Markov policies where each allocation
decision is a function of the state, at the arrival time of the item,
of the neighborhood where the item arrives. We will make this more
precise below.

Let $|S_i|=\kappa_i$ denote
the size of neighborhood $i$ and suppose the nodes in $S_i$ are
enumerated in some way, so that $S_i=\{s^{i}_1, \ldots, s^{i}_{\kappa_i}\}$.
\begin{df}
\label{connected}
We say that $j, k\in \{1, \ldots, n\}$ are neighbors  (and write
$j\sim k$), if $j, \,k\in S_i$ for some $i$. 
\end{df}
This equivalence relation can be used to define the graph ${\cal G}$
with vertices $\{1, \ldots, n\}$ and edges ${\cal E}$,
where $w=\langle j, k \rangle \in {\cal E}$ iff $j\sim k$.
Our main result (Theorem~\ref{time}) needs the following assumption.
\begin{cond}
\label{condS}
The graph ${\cal G}$ is  connected.
\end{cond}

Denote the configuration of the system at moment $t$ by
\[ X(t)=\bigl(X_1(t), \ldots, X_n(t)\bigr), \] 
where $X_i(t)$ is the number of items stored at node $i$ at time $t$.
The center of mass or average load of the configuration is
\[
M(t)=\frac{1}{n}\sum_{i=1}^n X_i(t),
\]
and we denote the \emph{shape} of the configuration by
\[
\tilde X(t)= \bigl(\tilde X_1(t), \ldots, \tilde X_n(t)\bigr) =
\bigl(X_1(t) - M(t), \ldots,  X_n(t) - M(t)\bigr),
\]  
the vector of loads relative to the center of mass. 
Note that, if a new item arrives at time $t$, then
$M(t)=M(t-)+\frac{1}{n}$. Also, if we know the shape $\tilde X(t)$, it
implies that we know which node is minimally loaded and we know the load differences between the nodes
(as $X_i(t)-X_j(t)=\tilde X_i(t)-\tilde X_j(t)$).

Obviously, the process $X(t)$ is Markovian for any decision rule that
depends only on the current node loads.  In order for the process
$\tilde X(t)$ to be Markovian, we require that the decision of
choosing the node is made accordingly to some decision rule which
depends only on the current {\it shape} of the system. Also, we are
mainly interested in {\em local} decision rules, that is, if an
item arrives to the set $S_i$, then {\it the only information
  about the configuration of the system that can be used to make a
  decision is what happens in the set $S_i$}. For example, the decision can be based on the 
differences $\tilde X_l(t)-\tilde X_j(t)$, $l, j\in S_i$.

If the decision rule is
  configuration independent and time homogeneous this gives rise
  to  a space homogeneous $(n-1)$-dimensional
random walk, which  is transient for $n>3$ and at best null recurrent for $n\le 3$.
 Therefore, if one wants positive recurrence in shape, the decision rule must depend 
on current configuration.  Of course, all nodes must receive arrivals for ergodicity
in shape to be achieved, hence the walk cannot live in a lower dimensional sub-space.
So, our goal is to find a rule for redistributing
 the arriving items at each moment of time  in a way to have 
positive recurrence in shape.
One of the possible  choices is to send the item to the node with minimal load
$S_i$ (Join the Shortest Queue routing policy).

We  present four routing policies. Two ensure the same rate of the
arrivals to different nodes, and the two others guarantee  stability in
shape, if some explicit conditions are fulfilled.
We note also that the conditions we refer to can be easily checked in
practice and the implementation of routing policies we propose is
algorithmically simple.

The paper is organized as follows. In  Section~\ref{notation}  we
introduce the notations and define the routing policies, in Section~\ref{results} we  state the
results.  In Section~\ref{prelim} we formulate the known facts we will
use in our proofs.  In Section~\ref{proofs1}, we  prove
 Theorem~\ref{time}, for which we need
some auxiliary lemmas, and then we prove Theorem~\ref{erg_alpha}.
In Section~\ref{s:noterg}, we first prove a lemma that translates the
condition of Theorem~\ref{t:noterg} into the language of convex
analysis, and prove Theorem~\ref{t:noterg}.

\section{Notations and definitions}
\label{notation}
Let us first introduce some notation. For $i=1,\dots, K$ denote
by $\Lambda_i$ the set of points 
$p^{(i)}=(p^{(i)}_1,\ldots, p^{(i)}_{\kappa_i}) \in \R^{\kappa_i}$ such that
\begin{equation}
\label{Lambda_i}
\begin{cases}
p^{(i)}_j\ge 0 & \text{for $j=1,\dotsc,\kappa_i$}, \\
\sum_{j=1}^{\kappa_i} p^{(i)}_j=1.
\end{cases}
\end{equation}

We use rather standard convention  that a vector $x\ge 0$ if its components are non-negative, and $x>0$ if its components are strictly positive.  


By $F$ denote the linear transformation that takes a point $x \in
\R^n$ to the point $y \in \R^n$ such that
$y_i=x_i-\frac{1}{n}\sum_{j=1}^n x_j$ for $i=1,\dots,n$.  In words,
the point $y$ represents deviations from the center of mass for the
configuration $x$. 
Let 
\begin{equation}
\M=F(\R^n)=\Bigl\{ y\in \R^n: \sum_{i=1}^n
y_i=0\Bigr\}. 
\end{equation}
%
%
Let $\tilde X(t)=F[X(t)]$.
The state space of the process $\tilde X(t)$ is  
\begin{equation}
{\cal M}=F(\N^n)=\Big\{y \in \big(n^{-1} \Z\big)^n: \quad \sum_{i=1}^n y_i=0\Big\}.
\end{equation}
Therefore, $\mathcal{M} \subset \M$.
We can say
informally that the dimension of the process $\tilde X(t)$ is $1$ less
than the dimension of $X(t)$.

A point $x=(x_1,\dotsc,x_n) \in \N^n$ represents the load of the system. By
$x_{S_i}$ denote the load of the nodes in $S_i$. Let $\1$ be the vector with all
ones: $\1=(1,\dotsc,1) \in \N^n$.

Now we define the notion of routing policy (RP).

\begin{df}
\label{routing_policy}
A {\em routing policy} $P$ is a function that takes a configuration $ x \in \N^n$ to
a point 
\(
P(x)=\bigl(P^{(1)}(x),\ldots, P^{(K )}(x)\bigr) \in \prod_{i=1}^K  \Lambda_i.
\)
For  the process $X(t)$  (or $\tilde X(t)$) with
routing policy~$P$,  an item 
arriving at neighbourhood $S_i$, when the configuration of the system is $x$,
is routed to  node $s^i_j$
with probability $P^{(i)}_j(x)$. The decisions are made 
  independently for each arrival. 
\end{df}

For the process $\tilde X(t)$ to  be Markovian, we suppose that all routing policies  satisfy the  following
\begin{cond}
The routing
policy~$P$ depends only on the current configuration shape, that is, for any admissible $x$ and $c\in \Z$ 
we have $P(x+c \1)=P(x)$.
\end{cond}


The decision about routing can be made using the complete information about configuration shape, or only partial information: 
\begin{df}
We say that a
routing policy $P$ is {\em local} if, for $i=1,\dots, K$,  the function
$P^{(i)}(x)$ depends only on the load of the nodes in $S_i$: 
for any $x$ and $y$ such that $x_{S_i}=y_{S_i}$, we have $P^{(i)}(x)=P^{(i)}(y)$.
\end{df}



In this paper we will consider four local routing policies.

\begin{df}
An {\em equilibrium routing policy (ERP)} is  a routing policy~$P$
such that~$P$ does not depend on $x$ and 
the resulting  arrivals at all  nodes
are independent Poisson processes with the same rate $1/n$
(recall that
$\sum_{i=1}^K  \lambda_i=1$).
\end{df}

\begin{df}
\label{df:serp}
A {\em strong equilibrium routing policy (SERP)} is 
an ERP with $P>0$.
\end{df}

Let us consider the following system of linear equations:
\begin{equation}
\label{lin_pr}
\begin{cases}
\sum_{j=1}^{\kappa_i} \alpha_{ij}=\lambda_i & \text{for $i=1,\dotsc, K$}, \\
\sum_{i=1}^K  \sum_{j=1}^{\kappa_i} \alpha_{ij} \delta_{\ell,s^i_j}=\frac{1}{n}
& \text{for $\ell=1,\dotsc,n$},
\end{cases}
\end{equation}
where $\delta_{\ell,m}$ is a Kronecker delta.

\begin{rmk}
 \label{exist_solution}
The system~\eqref{lin_pr} is a  special case of the maximum bipartite matching problem  and necessary and sufficient   conditions for existence of positive/non-negative solutions are well-known.
 
  For each non-empty collection of neighbourhoods $J \subset \{1, 2,
  \ldots, K  \}$ let $S_J = \cup_{j \in J} S_j$ and let $n_J$
  denote the number of nodes in $S_J$. Then, 
\begin{equation} 
\label{cond_lin_pr}
 \sum_{j \in J} \lambda_j \leq n_J/n \qquad\text{for all } J\subset\{1, \ldots, K  \}
\end{equation}
  is  necessary and sufficient for existence of non-negative solutions to~\eqref{lin_pr}. 
Strict inequality in \eqref{cond_lin_pr} for all $J$ except $\emptyset$ and $\{1, 2,
\ldots K \}$ is necessary and sufficient for the existence of
positive solutions to~\eqref{lin_pr}.

Indeed, if~\eqref{cond_lin_pr} is not satisfied, then at least one node in some $S_J$ must receive items at rate  greater than $1/n$, under any routing policy. The sufficiency can be shown using maximum-flow minimum-cut method (cf., for example,~\cite{FF1962, R1984}).

\end{rmk}

\begin{rmk}
\label{redirect} Note that
for any parameters of the model $S_1, \ldots, S_K $ and $\lambda_1, \dots, \lambda_n$ we have:
\begin{itemize}
 \item there exists an ERP iff~\eqref{lin_pr} has a non-negative solution;
\item there exists a SERP iff~\eqref{lin_pr} has a positive solution.
\end{itemize}
Indeed, if~\eqref{lin_pr} has a non-negative/positive solution we can define 
\[
P^{(i)}_j(x)=\alpha_{ij}/\lambda_i.\]
If we have an ERP/SERP, then
\[
\alpha_{ij}=\lambda_i P^{(i)}_j(x)
\]
is a non-negative/positive solution of~\eqref{lin_pr}. 
\end{rmk}

 We also rewrite this statement in the language of convex analysis (see Lemma~\ref{lm:riD}).

As solving~\eqref{lin_pr} is a problem of linear programming,  the existence of SERP can be easily checked in practice.

\begin{ex}
\noindent 
\begin{itemize} 
 \item 
Consider  a system with $n=3$ nodes and all possible neighborhoods of size~$2$, 
$\lambda_1+\lambda_2+\lambda_3=1$.
Then, there exists a positive solution of the system \eqref{lin_pr}  iff $\lambda_i<2/3$ for $i=1,2,3$.
\item Similarly, for  $n = 4$ and all possible  neighbourhoods
  of size~$2$, there exists a positive solution of the system \eqref{lin_pr}  
iff $\lambda_j < 1/2$, $j=1, \ldots, 6$, and
  $\sum_{j \in J} \lambda_j <  3/4$ for all $J$ such that $n_J = 3$. 
\end{itemize}
\end{ex}

Now we define the  other two routing policies which we study. 
For $x\in\N^n$, let
\begin{equation}
\label{def-s-max}
s^i_{j_{\max}}(x)=\max\Big\{s^i_{j}\in S_i: \;
  x_{s^i_j}=\max_{l=1, \ldots, \kappa_i}\{x_{s^i_l}\}\Big\}
\end{equation}
and 
\begin{equation}
\label{def-s-min}
s^i_{j_{\min}}(x)=\min\Big\{s^i_{j}\in S_i: \;
  x_{s^i_j}=\min_{l=1, \ldots, \kappa_i}\{x_{s^i_l}\}\Big\}.
\end{equation}
In words, for any load of the system $x\in \N^n$,  $s^i_{j_{\min}}(x)$ is the first node in $S_i$ such that in this node the load is minimal,
 $s^i_{j_{\max}}(x)$ is the last node in $S_i$ such that in this node the load is maximal.

\begin{df}
{\em Join the Shortest Queue (JSQ) routing policy} is
the routing policy~$P(x)=\bigl(P^{(1)}(x),\dotsc,P^{(K)}(x)\bigr)$, where
\begin{align*}
 P^{(i)}_j(x)=\left\{\begin{array}{ll} 
1 &\text{ if } s^i_{j}=s^i_{j_{\min}}(x),\\
0 &\text{otherwise}.
\end{array}
\right.
\end{align*}
 
\end{df}



\begin{df} 
Suppose that there exists a positive solution
    $\alpha_{ij}$ of~\eqref{lin_pr}. 
Let $0<\eps<\min{\alpha_{ij}}$. 
We define {\em $\eps$-perturbed strong equilibrium routing policy
    ($\eps$-PSERP)} 
as $P(x)=\bigl(P^{(1)}(x),\dotsc,P^{(K)}(x)\bigr)$, where
\begin{align*}
 P^{(i)}_j(x)=\left\{\begin{array}{ll} 
\frac{\alpha_{ij}+\eps}{\lambda_i} &\text{ if } s^i_{j}=s^i_{j_{\min}}(x),\\
\frac{\alpha_{ij}-\eps}{\lambda_i} &\text{ if } s^i_{j}=s^i_{j_{\max}}(x),\\
\frac{\alpha_{ij}}{\lambda_i} &\text{otherwise}.
\end{array}
\right.
\end{align*}
If $\kappa_i=1$ (i.e., the neighborhood $S_i$ has size $1$), then we have no freedom to choose probabilities and
 $P^{(i)}_j(x)=1$ for any $x$.
\end{df}

Note that 
in each of the four cases the routing policy can be chosen to be local.
Indeed, in the case of JSQ it is clear immediately from the definition.
In each of the other three cases, we first need to note that  we can choose the same
solution of \eqref{lin_pr} for all $x \in \N^n$, then it is easy to see that the
corresponding policy is local.
Moreover, in the cases of ERP and SERP it does not depend on $x$.

\medskip


We study the behavior of the process $\tilde X(t)$ that has state space $\mathcal{M}$. 
In order to simplify the notation, we prefer to keep the same symbol
for the process with any  RP; instead when dealing with $X(t)$ or $\tilde X(t)$
we will state explicitly which RP is used.

Let $\{X^e(m)\}_{m\in \N}$ (resp. $\{\tilde X^e(m)\}_{m\in \N}$) be the embedded Markov chain for the process $\{X(t)\}_{t\ge 0}$ 
(resp. $\{\tilde X(t)\}_{t\ge 0}$), obtained 
when we look at the system only at the moments of arrivals. Note that $\{X^e(m)\}_{m\in \N}$  and 
$\{\tilde X^e(m)\}_{m\in \N}$ are indeed Markov chains,
as the arrivals are Poisson. 
Note also that $\{\tilde X^e(m)\}_{m\in \N}$ has period $n$ under any of the policies considered (indeed, if  $\tilde X^e(m)=\tilde x$, we need the same number of items to arrive at every node to obtain  $\tilde X^e(m')=\tilde x$, so we must have $m'= n l$ for some $l$). 
For ERP, SERP and $\eps$-PSERP the process $\{\tilde X^e(m)\}_{m\in \N}$ is irreducible, as
all nodes have positive arrival rates and thus any shape can be obtained from any other shape. 
The situation is more delicate for JSQ routing policy. For example, with JSQ, if node $j$ does not belong to a neighborhood of size $1$, then starting from configuration $\tilde X^e(0)=0$
it is impossible to obtain configuration with  $\tilde X_i^e(m)=\tilde x$ for all $i\ne j$ and 
$\tilde X_j^e(m)=\tilde x+2/n$. It important to note, however, that the configuration  $\tilde X^e(m)=0$ {\it is reachable from any configuration}.  

 By $\tau$ denote the time of the first return to the origin:
\begin{equation}
\label{tau}
\tau=\inf\{m>0: \tilde X^e(m)=0\}.
\end{equation}
We say that
\begin{itemize}
\item[(a)]
$\{\tilde X^e(m)\}_{m\in \N}$ is transient if $\P(\tau=\infty \mid \tilde X^e(0)=0)>0$,
\item[(b)]
 $\{\tilde X^e(m)\}_{m\in \N}$ is recurrent if $\P(\tau<\infty \mid  \tilde X^e(0)=\tilde x)=1$ 
for any $\tilde x\in {\cal M}$, 
\item[(c)]
 $\{\tilde X^e(m)\}_{m\in \N}$ is positive recurrent if $\E (\tau\mid\tilde X^e(0)=\tilde x ) <\infty$
 for any $\tilde x\in {\cal M}$.
\end{itemize}
We prefer to give the definition in this form
because, as we will see below, (b) and (c) either hold for all or for no $\tilde x \in  {\cal M}$. 
%
\section{Recurrence/transience classification}
\label{results}
Since the rates of our processes are bounded away from $0$ and $\infty$,
positive recurrence  of $\{\tilde X(t)\}_{t\ge 0}$  is equivalent to positive recurrence
of $\{\tilde X^e(m)\}_{m\in \N}$. So, we will prove the results for
$\{\tilde X^e(m)\}_{m\in \N}$.

Define the shape magnitude as
\begin{equation}
\label{def_disp}
{\mathfrak D}\bigl(\tilde X(t)\bigr)=\sum_{i=1}^n \bigl(\tilde X_i(t)\bigr)^2=\sum_{i=1}^n \bigl(X_i(t)-M(t)\bigr)^2
\end{equation}
$\big(\text{so ${\mathfrak D}\bigl(\tilde X(t)\bigr)$ is in fact the square of the Euclidean norm of 
$\tilde X(t)$}\big)$.

\begin{theo}
\label{time}
Suppose that Condition~\ref{condS} is satisfied
and there exists a positive solution of~\eqref{lin_pr}. 
\begin{itemize}
 \item[(i)]
Suppose that we  construct the process $\{\tilde X^e(m)\}_{m\in \N}$ using either
JSQ routing policy or $\eps$-PSERP.
Then  $\tilde X^e(m)$ is positive recurrent. Moreover, there exists $c>0$ such
that for all $0<c'<c$ we have 
\[
\E (e^{c'\tau}\mid\tilde X^e(0)=x)<\infty
\]
 for all $x$.
\item[(ii)]  Also, JSQ routing policy minimizes the expected shape magnitude, that is, 
for any routing policy we have
\begin{eqnarray*}
\lefteqn{\E^{\mbox{\tiny any RP}} \bigl[ {\mathfrak D}\bigl(\tilde X^e(m+1)\bigr)\mid \tilde X^e(m)=x\bigr]}\\
&\ge& \E^{\mbox{\tiny JSQ}} \bigl[ {\mathfrak D}\bigl(\tilde X^e(m+1)\bigr)\mid \tilde X^e(m)=x\bigr].
\end{eqnarray*}
\end{itemize}
\end{theo}

Note that using ERP or SERP it is impossible to have positive recurrence of  $\tilde X^e(m)$. Indeed, these routing policies provide independent Poisson arrivals with the same rate to all nodes. Then  the behavior of the shape can be described by a $(n-1)$-dimensional random walk with zero drift, which is transient if $n>3$ and null-recurrent if $n\le 3$.

If the Condition~\ref{condS} is not fulfilled, then we have two or more disconnected components, that is, sets of nodes such that arrivals to one of these sets cannot be routed to the other. In this case, it is impossible to obtain positive recurrence in shape, for any routing policy. If the number of disconnected components is at least 4, then even null-recurrence is impossible (as in the argument above).



We also have the following converse results (in some sense) to
Theorem~\ref{time}.  Note that in
Theorems~\ref{erg_alpha} and~\ref{t:noterg} we do not require the
routing policy $P$ to be local.

\begin{theo}
\label{erg_alpha} Fix the parameters of the model: $S_1, \ldots, S_K $, 
$\lambda_1, \ldots, \lambda_K $.
Suppose that there exists a routing policy $P$ such that the process $\tilde X(t)$ with the routing policy $P$ is recurrent. 
Then there exists a non-negative solution of~\eqref{lin_pr} (and thus
for the model with these parameters there exists an ERP).
\end{theo}

We can also rewrite Theorem~\ref{erg_alpha} in a different way:
\begin{cor}
\label{transient}
Fix the parameters of the model: $S_1, \ldots, S_K $, 
$\lambda_1, \ldots, \lambda_K $.
Suppose that there is no non-negative solution $\alpha_{ij}$ of the system 
\eqref{lin_pr}. Then for any routing policy $P$, the process $\tilde X(t)$
with the routing policy $P$ is transient.  
\end{cor}

\begin{theo}
\label{t:noterg}
Fix the parameters of the model: $S_1, \ldots, S_K $, 
$\lambda_1, \ldots, \lambda_K $.
Suppose that there is no positive solution $\alpha_{ij}$ of the system
\eqref{lin_pr}.
Then for any routing policy $P$, the process $\tilde X(t)$ with the routing policy $P$ is not positive recurrent.
\end{theo}

The following problem is still open. Fix the parameters of the model: $S_1,
\ldots, S_K $, $\lambda_1, \ldots, \lambda_K $.
Suppose that there is no positive solution $\alpha_{ij}$ of the system 
\eqref{lin_pr}, but  there exists a non-negative solution. Under which
conditions on the parameters of the model $S_1, \ldots, S_{\cal
K}$, $\lambda_1, \ldots, \lambda_K $ (and $n$) does there exist a (local) routing
policy $P$ such that the process $\tilde X(t)$ with the routing policy $P$ is recurrent?

\section{Proofs}
The structure of this section is as follows. First (Section~\ref{prelim})  we formulate some known fact which we will use in our proofs. In Section~\ref{proofs1}, we introduce some notations and define two functions ($f$ and $g$) we will use to prove~Theorem~\ref{time}. Then we prove four lemmas, obtaining bounds on 
\[
\E\bigl[f\bigl(X^e(m+1)\bigr)-f\bigl(X^e(m)\bigr)\mid X^e(m)=x\bigr]
\]
 for JSQ and $\eps$-PSERP. Using these bounds, we prove  Theorem~\ref{time}. Then we prove Theorem~\ref{erg_alpha}.
In Section~\ref{s:noterg}, we first recall some definitions from complex analysis and apply these to our model. Then we prove Lemma~\ref{lm:riD}, which  translates the
condition of Theorem~\ref{t:noterg} into the language of convex
analysis, and then we finish the proof of Theorem~\ref{t:noterg}.

\subsection{Preliminaries}
\label{prelim}
We state some known results that we will use in our proofs. Note that
Theorems~\ref{t:ergodic} and~\ref{t:non-ergodic}  are Theorems~2.2.3 and 2.2.6
respectively from~\cite{FMM}, where we use `positive recurrent' instead of
`ergodic'. This change is necessary as our Markov chains are periodic. That the results also hold 
for periodic chains is mentioned in Section~1.1 of~\cite{FMM}. 
In fact, to see that the reformulated theorems are valid it suffices to
consider the Markov chain $\eta_{\ell}$ at embedded instants $\ell=k+pr$, where $p$ is
the period of the chain and $k$ is a fixed number.

Let us consider a time homogeneous irreducible  Markov chain
$\eta_m$ with
  countable state
space $\cal H$.
\begin{theo}
\label{t:ergodic}
 The Markov chain $\eta_m$ is positive recurrent if and
only if there exists a positive function $f(x), x \in \cal H$, a number
$\varepsilon
> 0 $ and a finite set $A \in \cal H$ such that for every $m$ we have
\begin{eqnarray}
\label{e:ergodic0}
{\bf E}[f(\eta_{m+1})-f(\eta_m) \mid \eta_m=x] &\le& -\varepsilon, \quad x\notin A,\\
{\bf E}[f(\eta_{m+1}) \mid \eta_m=x] &<&\infty, \quad x\in A. \nonumber
\end{eqnarray}
\end{theo}

\begin{theo}
\label{t:non-ergodic}
For the Markov chain $\eta_m$ to be not positive recurrent,
it is sufficient that there exists a function $f(x), x
\in \cal H$, and constants $C\in \R$ and $d>0$ such that
\begin{itemize}
 \item
  for every $m$ we have \[{\bf E}[f(\eta_{m+1})-f(\eta_m)
  \mid \eta_m=x] \ge 0, \; x \in \{f(x)>C\},\]
     where the sets $ \{x\mid f(x)>C\}$ and $ \{x\mid
f(x)\le C\}$ are non empty;
 \item for every $m$ we have
 \[{\bf E}\bigl[|f(\eta_{m+1})-f(\eta_m)| \mid \eta_m=x\bigr]\le d, \; x\in \cal H.\]
\end{itemize}
\end{theo}

The following theorem is an immediate consequence of Theorem 2.1.7 from~\cite{FMM}.
\begin{theo}
\label{217}
Let $(\Omega, {\cal F}, \P)$ be the probability space and
$\{{\cal F}_n, \; n\ge 0\}$ be an increasing family of
$\sigma$-algebras. Let $\{{\mathfrak S}_l, \; l\ge 0\}$ be a sequence of random
variables such that ${\mathfrak S}_l$ is  ${\cal F}_l$-measurable, and
${\mathfrak S}_0$ is a constant. Let
\[y_{k+1}={\mathfrak S}_{k+1}-{\mathfrak S}_k.
\]
If there exist
positive numbers
$\eps$, $M$, such that for each $k$ we have
\[\E[y_{k+1}\mid {\cal F}_k]\le -\eps, \mbox{ a.s.}\]
\[|y_{k+1}|<M \mbox{ a.s.},\]
then, for any $\delta_1<\eps$, there exist constants
$C=C({\mathfrak S}_0)$ and $\delta_2>0$, such that,
for any $m>0$,
\[\P[{\mathfrak S}_m>-\delta_1 m]<Ce^{-\delta_2 m}.\]
\end{theo}


\subsection{Proofs of Theorems~\ref{time} and \ref{erg_alpha}}
\label{proofs1}
To prove Theorem~\ref{time},
we need some additional notations and four lemmas.  

Suppose that we are using either
JSQ   routing policy or
$\eps$-PSERP to
construct the process $X^e(m)$ (for now, it does not matter which one).
We are going to construct a supermartingale with bounded jumps,  
that will allow us to obtain exponential bounds on $\tau$ (see~(\ref{tau})
for the definition of $\tau$)
and thus to prove positive recurrence of  $\tilde X^e(m)$.

Let \[f\bigl(X^e(m)\bigr)=f\bigl(X_1^e(m), \ldots, X_n^e(m)\bigr)=\sum_{i=1}^n \bigl(X_i^e(m)-M^e(m)\bigr)^2={\mathfrak D}\bigl(\tilde X^e(m)\bigr),\]
where ${\mathfrak D}\bigl(\tilde X^e(m)\bigr)$ is the shape magnitude defined in~(\ref{def_disp})
and
\[g\bigl(\tilde X^e(m)\bigr)=\sqrt{f\bigl(X^e(m)\bigr)}=\Bigl(\sum_{i=1}^n
\bigl(X_i^e(m)-M^e(m)\bigr)^2\Bigr)^{1/2}.\]
We will prove that $g\bigl(\tilde X^e(m)\bigr)$ is a supermartingale with bounded jumps. 
To do that, we will need some auxiliary lemmas.
In Lemmas~\ref{lemmaPERP} and~\ref{lemmaJSQ} we estimate 
$ \E\bigl[f\bigl(X^e(m+1)\bigr)-f\bigl(X^e(m)\bigr)\mid X^e(m)=x \bigr]$ in terms of $X^e(m)$ for
$\eps$-PSERP and JSQ respectively. In Lemma~\ref{lemma_g1} we obtain a bound on
 $\bigl|f\bigl(X^e(m+1)\bigr)-f\bigl(X^e(m)\bigr)\bigr|$, which is needed for the proof that
$g\bigl(\tilde X^e(m)\bigr)$ has bounded jumps.

First, we introduce
the process $\bigl(Y_1(m), \ldots, Y_n(m)\bigr)$  obtained when the
item that arrives at $S_i$
is directed to  node $s^i_j$ with probability
$p_j^{(i)}=\alpha_{ij}/\lambda_i$, $j=1, \ldots \kappa_i$ (that is, using SERP).
The processes
$X^e(m)$ and $Y(m)$ are defined in  the same probability space, 
use the same  arrivals, and if $X^e(m)=Y(m)=x $,
then $X^e(m+1)$ and $Y(m+1)$ are obtained from $x $ using the respective
routing policies (independently for $X^e(m+1)$ and $Y(m+1)$). In addition,  it is clear that 
$\P\bigl(Y(m)=x\bigr) > 0$ iff $ \P\bigl(X^e(m)=x\bigr)>0$. 

Using the fact that $\alpha_{ij}$'s are such that arriving  items are routed to  node $i$
with probability $1/n$ for any $i$,
we have
\begin{eqnarray}
\lefteqn{\E\bigl[\bigl(Y_i(m+1)-M^Y(m+1)\bigr)^2-\bigl(Y_i(m)-M^Y(m)\bigr)^2\mid Y(m)\bigr]}\nonumber\\
&=&\frac{1}{n}\Big(\Big(Y_i(m)+1-M^Y(m)-\frac{1}{n}\Big)^2-\bigl(Y_i(m)-M^Y(m)\bigr)^2\Big)\nonumber\\
&&{}+\frac{n-1}{n}
\Big(\Big(Y_i(m)-M^Y(m)-\frac{1}{n}\Big)^2-\bigl(Y_i(m)-M^Y(m)\bigr)^2\Big)\nonumber\\
&=&\frac{1}{n}-\frac{1}{n^2}\label{oc_Yi},
\end{eqnarray}
where $M^Y(m)=\frac{1}{n}\sum_{k=1}^n Y(k)$, as $M^Y(m+1)=M^Y(m)+\frac{1}{n}$.
Thus,
\begin{equation}
\label{oc_Y}
\E\bigl[f\bigl(Y(m+1)\bigr)-f\bigl(Y(m)\bigr)\mid Y(m)\bigr]=n\Big(\frac{1}{n}-\frac{1}{n^2}\Big)=1-\frac{1}{n}.
\end{equation}

Denote by $C_i$ the event that an item arrives at set $S_i$.
Recall~\eqref{def-s-max} and~\eqref{def-s-min}. From now on, in order to simplify notation, instead of writing
$s^i_{j_{\max}}\bigl(X^e(m)\bigr)$ and $s^i_{j_{\min}}\bigl(X^e(m)\bigr)$, we
will write  $s^i_{j_{\max}}$ and $s^i_{j_{\min}}$.
Also, instead of $X_{s^i_j}^e(m)$ we will write $X^e(i,j,m)$; analogously 
for $\tilde X^e(t)$, $Y(m)$  and $\tilde Y(m)$.

\begin{lm}
\label{lemmaPERP}
Suppose that the process  $\{X^e(m)\}_{m\in \N}$ is constructed using $\eps$-PSERP.
Then
\begin{eqnarray}
\lefteqn{\E\bigl[f\bigl(X^e(m+1)\bigr)-f\bigl(X^e(m)\bigr)\mid X^e(m)=x \bigr]}\nonumber\\
&=& -2\eps\sum_{i=1}^K  \bigl(X^e(i,j_{\max},m)-X^e(i,j_{\min},m)\bigr)+
 \label{oc2_X} 1-\frac{1}{n}.
\end{eqnarray}
\end{lm}

\medskip
\noindent
{\it Proof of Lemma~\ref{lemmaPERP}.}
Suppose $|S_i|>1$.
We have, for $x$ such that $\P\bigl(X^e(m)=x\bigr)>0$ $\big(\text{and thus $\P\bigl(Y(m)=x\bigr) > 0$}\bigr)$,
\begin{eqnarray}
\lefteqn{\E\bigl[f\bigl(X^e(m+1)\bigr)-f\bigl(Y(m+1)\bigr)\mid X^e(m)=Y(m)=x , \; C_i\bigr]}\nonumber\\
&=& \E\Big[ \sum_{j=1}^{\kappa_i}\bigl(\tilde X^e(i,j,m+1)\bigr)^2-\bigl(\tilde Y(i,j,m+1)\bigr)^2
            \;\bigl|\;  X^e(m)=Y(m)=x , \; C_i\Big]\nonumber\\
&=&\frac{\eps}{\lambda_i}  \biggl(\Big(X^e(i,j_{\min},m)+1-M^e(m)-\frac{1}{n}\Big)^2\nonumber\\
&& ~~~~~~~~ +\sum_{j\ne j_{\min}}\Big(X^e(i,j,m)-M^e(m)-\frac{1}{n}\Big)^2 \biggr)\nonumber\\
&&-\frac{\eps}{\lambda_i}  \biggl(\Big(X^e(i,j_{\max},m)+1-M^e(m)-\frac{1}{n}\Big)^2\nonumber\\
&&~~~~~~~~+\sum_{j\ne j_{\min}}\Big(X^e(i,j,m)-M^e(m)-\frac{1}{n}\Big)^2 \biggr)\nonumber\\
&=&-\frac{2\eps}{\lambda_i} \bigl(X^e(i,j_{\max},m)-X^e(i,j_{\min},m) \bigr) \label{oc_X}
\end{eqnarray}
as we conditioned  on $X^e(m)=Y(m)=x $.
Thus,
\begin{eqnarray}
\lefteqn{\E\bigl[f\bigl(X^e(m+1)\bigr)-f\bigl(X^e(m)\bigr)\mid X^e(m)=x , \; C_i\bigr]}\nonumber\\
&=&\E\bigl[f\bigl(X^e(m+1)\bigr)-f\bigl(Y(m+1)\bigr)\mid X^e(m)=Y(m)=x , \; C_i\bigr]\nonumber\\
&&+\E\bigl[f\bigl(Y(m+1)\bigr)-f\bigl(X^e(m)\bigr)\mid X^e(m)=Y(m)=x , \; C_i\bigr]\label{oc1_X}\\
&=& -\frac{2\eps}{\lambda_i} \bigl(X^e(i,j_{\max},m)-X^e(i,j_{\min},m) \bigr)\nonumber\\
&& +\E\bigl[f\bigl(Y(m+1)\bigr)-f\bigl(Y(m)\bigr)\mid Y(m)=x , \; C_i\bigr]\nonumber
\end{eqnarray}
and
\begin{eqnarray}
\lefteqn{\E\bigl[f\bigl(X^e(m+1)\bigr)-f\bigl(X^e(m)\bigr)\mid X^e(m)=x \bigr]}\nonumber\\
&=&\sum_{i=1}^K  \lambda_i \E\bigl[f\bigl(X^e(m+1)\bigr)-f\bigl(X^e(m)\bigr)\mid X^e(m)=x , \; C_i\bigr]\nonumber\\
&=& -2\eps\sum_{i=1}^K  \bigl(X^e(i,j_{\max},m)-X^e(i,j_{\min},m)\bigr)+
 \label{oc2_Xa} 1-\frac{1}{n}.
\end{eqnarray}
Note that if there is a neighborhood of size 1, by Condition~\ref{condS} it should be subset of another neighborhood, of size at least 2. As the terms corresponding to neighborhoods of size 1 in~\eqref{oc2_Xa} will be equal to 0, the equation~\eqref{oc2_Xa} still holds. 
Lemma~\ref{lemmaPERP} is proved.
\qed

\begin{lm}
\label{lemmaJSQ}
Suppose that the process $\{X^e(m)\}_{m\in \N}$ is constructed using JSQ routing policy.
Then
\begin{eqnarray}
\lefteqn{\E\bigl[f\bigl(X^e(m+1)\bigr)-f\bigl(X^e(m)\bigr)\mid X^e(m)=x \bigr]}\nonumber\\
&=& -2\sum_{i=1}^K  \sum_{j\ne j_{\min}}\alpha_{i j} \bigl(X^e(i,j,m)-X^e(i,j_{\min},m)\bigr)
 +1-\frac{1}{n}\label{oc2a_X}.
\end{eqnarray}
\end{lm}

\medskip
\noindent
{\it  Proof of Lemma~\ref{lemmaJSQ}}. 
Analogously to~(\ref{oc_X}),
\begin{eqnarray}
\lefteqn{\E\bigl[f\bigl(X^e(m+1)\bigr)-f\bigl(Y(m+1)\bigr)\mid X^e(m)=Y(m)=x , \; C_i\bigr]}\nonumber\\
&=&\sum_{j\ne j_{\min}}\frac{\alpha_{i j}}{\lambda_i}\biggl(\Big(X^e(i,j_{\min},m)+1-M^e(m)-\frac{1}{n}\Big)^2\nonumber\\
&&+\sum_{j''\ne j_{\min}}\Big(X^e(i,j'',m)-M^e(m)-\frac{1}{n}\Big)^2\nonumber \\
&&-\Big(\Big(Y(i,j,m)+1-M^e(m)-\frac{1}{n}\Big)^2 +
     \sum_{j'\ne j}\Big(Y(i,j',m)-M^e(m)-\frac{1}{n}\Big)^2 \Big)\biggr)\nonumber
\\
&=& \sum_{j\ne j_{\min}}\frac{\alpha_{i j}}{\lambda_i}
      \biggl(\Big(X^e(i,j_{\min},m)+1-M^e(m)-\frac{1}{n}\Big)^2\nonumber\\
&&+ \Big(X^e(i,j,m)-M^e(m)-\frac{1}{n}\Big)^2 \nonumber\\
&&-\Big(\Big(Y(i,j,m)+1-M^e(m)-\frac{1}{n}\Big)^2 +
 \Big(Y(i,j_{\min},m)-M^e(m)-\frac{1}{n}\Big)^2 \Big)\biggr)\nonumber\\
&=& - \sum_{j\ne j_{\min}}\frac{2\alpha_{i j}}{\lambda_i}\bigl( X^e(i,j,m)-X^e(i,j_{\min},m)\bigr)\label{oc3a_X}.
\end{eqnarray}
So,
\begin{eqnarray}
\lefteqn{\E\bigl[f\bigl(X^e(m+1)\bigr)-f\bigl(X^e(m)\bigr)\mid X^e(m)=x , \; C_i\bigr]}\nonumber\\
&=& - \sum_{j\ne j_{\min}}\frac{2\alpha_{i j}}{\lambda_i}\bigl( X^e(i,j,m)-X^e(i,j_{\min},m)\bigr)\nonumber\\
&&+\E\bigl[f\bigl(Y(m+1)\bigr)-f\bigl(Y(m)\bigr)\mid Y(m)=x , \; C_i\bigr]\label{oc1a_X}
\end{eqnarray}
and
\begin{eqnarray}
\lefteqn{\E\bigl[f\bigl(X^e(m+1)\bigr)-f\bigl(X^e(m)\bigr)\mid X^e(m)=x \bigr]}\nonumber\\
&=&\sum_{i=1}^K  \lambda_i \E\bigl[f\bigl(X^e(m+1)\bigr)-f\bigl(X^e(m)\bigr)\mid X^e(m)=x , \; C_i\bigr]\nonumber\\
&=& -2\sum_{i=1}^K  \sum_{j\ne j_{\min}}\alpha_{i j} \bigl(X^e(i,j,m)-X^e(i,j_{\min},m)\bigr)
 +1-\frac{1}{n}\label{oc2a_Xb}.
\end{eqnarray}
Lemma~\ref{lemmaJSQ} is proved.
\qed

\medskip
 Denote by $e_i$ the $i$-th coordinate vector, $i=1, \ldots, n$.
The next lemma will be used to bound
  jumps in $f$ due to any possible one-step changes to $x$. 

\begin{lm}
\label{lemma_g1}
Let $x\in \N^n$ and $m(x)=\frac{1}{n}\sum_{j=1}^n x_j$. If
$\sum_{i=1}^n \big(x_i-m(x)\big)^2>0$, then  for each $e_i$, $i=1, \ldots, n$, 
\begin{equation}
\label{eq1g0}
| f(x+e_i)-f(x) |\le 4 \sqrt{f(x)}.
\end{equation}
\end{lm}

\medskip
\noindent
{\it Proof of Lemma~\ref{lemma_g1}.}
Without loss of generality, consider the first coordinate vector $e_1$.
We have then
\begin{eqnarray}
\lefteqn{f(x+e_1)}\nonumber\\
&=&\Big(x_1+1-m(x)-\frac{1}{n}\Big)^2+\sum_{i=2}^n \Big(x_i-m(x)-\frac{1}{n}\Big)^2\nonumber\\
&=&\big(x_1-m(x)\big)^2+2\big(x_1-m(x)\big)
\Big(1-\frac{1}{n}\Big)+\Big(1-\frac{1}{n}\Big)^2\nonumber\\
&&+
\sum_{i=2}^n \big(x_i-m(x)\big)^2-\frac{2}{n}\sum_{i=2}^n \big(x_i-m(x)\big)+\frac{n-1}{n^2}\nonumber\\
&=&\sum_{i=1}^n \big(x_i-m(x)\big)^2+2\big(x_1-m(x)\big)+1-\frac{1}{n}\label{kvadr}
\end{eqnarray}
as
\[\frac{1}{n}\sum_{i=1}^n \big(x_i-m(x)\big)=0.\]
Hence for each $i=1, \ldots, n$ we have
\begin{equation}
\label{eq1f}
f(x+e_i)-f(x)=2\big(x_i-m(x)\big)+1-\frac{1}{n}.
\end{equation}
It remains to show that, if $\sum_{i=1}^n \big(x_i-m(x)\big)^2>0$, then
\[
 \Big|2\big(x_i-m(x)\big)+1-\frac{1}{n}\Big|\le 4 \sqrt{f(x)}.
\]
Note that  $\sum_{i=1}^n \big(x_i-m(x)\big)^2>0$ implies that there exists  $i,j$
such that $|x_i-x_j|\ge 1$. Thus, there exists at least one $l$ such that $|x_l-m(x)|\ge 1/2$ which implies that 
$\sqrt{f(x)}\ge 1/2$.
So,
\begin{align*}
 |f(x+e_i)-f(x)|&\le 2| x_i-m(x)|+1\\
&\le 2\Big[\sum_{i=1}^n \big(x_i-m(x)\big)^2\Big]^{1/2}+1\\
&\le 4\sqrt{f(x)}.
\end{align*}

Lemma~\ref{lemma_g1} is proved.
\qed

Lemma~\ref{lemma_g1} implies that, for any RP, if $\sum_{i=1}^n \big(X_i^e(m)-M^e(m)\big)^2>0$, then 
\begin{equation}
 \label{eq1g}
 \bigl| f\bigl(X^e(m+1)\bigr)-f\bigl(X^e(m)\bigr) \bigr|\le 4\sqrt{f\bigl(X^e(m)\bigr)} = 4 g\bigl(\tilde X^e(m)\bigr).
 \end{equation}

It is important to note that the next computations are valid for JSQ 
and for $\eps$-PSERP. 

\begin{lm}
\label{lemma_g2}
There exist  $c_2>0$ and $a>0$, such that for all $x\in \N^n$ with
 \[
\max_{i=1, \ldots, n}\Big|x_i-\frac{1}{n}\sum_{j=1}^n x_j\Big|\ge a
\] 
it holds that
\begin{equation}
\label{eqf2}
\E\bigl[f\bigl(X^e(m+1)\bigr)-f\bigl(X^e(m)\bigr) \mid X^e(m)=x\bigr] \le - c_2 \sqrt{f(x)}.
\end{equation}
\end{lm}

\medskip
\noindent
{\it Proof of Lemma~\ref{lemma_g2}.}
We have
\begin{eqnarray*}
f\bigl(X^e(m)\bigr)&=& \sum_{l=1}^n \big(X_l^e(m)-M^e(m)\big)^2\\
&\le& \sum_{i=1}^K \sum_{j\in S_i} \bigl(X_j^e(m)-M^e(m)\bigr)^2\\
&\le &  \sum_{i=1}^K |S_i| \;\max_{j\in S_i} \bigl\{\bigl(X_j^e(m)-M^e(m)\bigr)^2\bigr\}\\
&\le & n \sum_{i=1}^K  \max_{j\in S_i}\bigl\{ \bigl(X_j^e(m)-M^e(m)\bigr)^2\bigr\}.\\
\end{eqnarray*}

We now show that under Condition~\ref{condS} we have
\begin{equation}
\label{conds1}
\sum_{i=1}^K  \max_{j\in S_i} \bigl(X_j^e(m)-M^e(m)\bigr)^2
\le c_3 \Big(\sum_{i=1}^K  \bigl(X^e(i,j_{\max},m)-X^e(i,j_{\min},m)\bigr) \Big)^2
\end{equation}
and also
\begin{equation}
\label{conds2}
 \sum_{i=1}^K  \max_{j\in S_i} \bigl(X_j^e(m)-M^e(m)\bigr)^2
\le c_4 \Big(\sum_{i=1}^K  \sum_{j\ne j_{\min}}\alpha_{i j} \bigl(X^e(i,j,m)-X^e(i,j_{\min},m)\bigr)\Big)^2.
\end{equation}
Let us consider (\ref{conds1}).
If \[
X^e(i,j_{\min},m)\le M^e(m)\le X^e(i,j_{\max},m),
\]
then,
obviously, 
\[
\bigl(X_j^e(m)-M^e(m)\bigr)^2\le \bigl(X^e(i,j_{\max},m)-X^e(i,j_{\min},m)\bigr)^2.
\]
Suppose that $M^e(m)< X^e(i,j_{\min},m)$ (the case
$M^e(m)> X^e(i,j_{\max},m)$ can be treated analogously).
Consider the sets of nodes $\{j : X^e_j(m) \leq
  M^e(m)\}$ and $\{j : X^e_j(m) > M^e(m)\}$. By Condition 1.1 some
  neighbourhood contains nodes from each of these sets and hence there
  exists $i^*$ 
such that
\[
X^e(i^*,j_{\min},m)\le M^e(m)\le X^e(i^*,j_{\max},m)
\]
and a sequence of neighbourhoods indexed by $i_0=i, i_1, \ldots, i_k=i^*$ such that
$S_{i_{\ell-1}}\cap S_{i_\ell}\ne \emptyset$, $\ell=1, \ldots, k$.
Thus,
\begin{equation}
\label{eq_summa}
\begin{split}
&\max_{j\in S_i} \bigl|X_j^e(m)-M^e(m)\bigr|\\ 
&\quad\le X^e(i,j_{\max},m)-M^e(m)\\
&\quad\le X^e(i,j_{\max},m)-X^e(i,j_{\min},m) + X^e(i,j_{\min},m)-M^e(m)\\
&\quad\le X^e(i,j_{\max},m)-X^e(i,j_{\min},m) + X^e(i_1,j_{\max},m)-M^e(m).
\end{split}
\end{equation}
The last inequality is due to the fact that, as $S_i \cap S_{i_1}\ne \emptyset$,
it holds
\begin{eqnarray*}
\lefteqn{X^e(i,j_{\min},m)= \min_{s^i_j\in S_i} X^e(i,j,m)
   \le  \min_{s^i_j\in S_i\cap S_{i_1} } X^e(i,j,m)}\\
&\le& \max_{s^i_j\in S_i\cap S_{i_1} } X^e(i,j,m)
\le \max_{s^{i_1}_j\in S_{i_1}}X^e(i_1,j,m) =X^e(i_1,j_{\max},m) 
\end{eqnarray*}
Continuing~\eqref{eq_summa}, we get
\[
\begin{split} 
&\max_{j\in S_i} \bigl|(X_j^e(m)-M^e(m)\bigr| \\
&\quad\le X^e(i,j_{\max},m)-M^e(m)\\
&\quad\le  X^e(i,j_{\max},m)-X^e(i,j_{\min},m) + X^e(i_1,j_{\max},m)-M^e(m) \\
&\quad\le  X^e(i,j_{\max},m)-X^e(i,j_{\min},m) 
+ X^e(i_1,j_{\max},m)- X^e(i_1,j_{\min},m)\\
&\qquad+ X^e(i_2,j_{\max},m)-M^e(m)
\end{split}
\]
and so on until $i_k=i^*$ (at the last step one has to use $X^e(i^*,j_{\min},m)\le M^e(m)$). 
So, we obtain
\begin{eqnarray*}
\max_{j\in S_i} \bigl|(X_j^e(m)-M^e(m)\bigr|
\le \sum_{\ell=0}^{k} \bigl(X^e(i_\ell,j_{\max},m)-X^e(i_\ell,j_{\min},m)\bigr)
\end{eqnarray*}
and  \eqref{conds1} follows with some $c_3\le K $. The argument for \eqref{conds2}  is similar.

Then  Lemma~\ref{lemmaPERP} together with (\ref{conds1}) (for $\eps$-PSERP), and Lemma~\ref{lemmaJSQ} together with (\ref{conds2}) (for JSQ),
imply that, for some $c_2>0$,
\begin{equation}
\label{eqf2a}
\E\bigl[f\bigl(X^e(m+1)\bigr)-f\bigl(X^e(m)\bigr)\mid X^e(m)\bigr] \le - c_2 \sqrt{f\bigl(X^e(m)\bigr)}=-c_2g\bigl(\tilde X^e(m)\bigr),
\end{equation}
when $f\bigl(X^e(m)\bigr)$ is large enough.
Lemma~\ref{lemma_g2} is proved. 
\qed

\medskip
{\it Proof of Theorem~\ref{time}}.
First, we verify that $g\bigl(\tilde X^e(m)\bigr)$ has bounded jumps.
If $\sum_{i=1}^n \big(X_i^e(m)-M^e(m)\big)^2=0$, then, obviously,
$g\bigl(\tilde X^e(m+1)\bigr)-g\bigl(\tilde X^e(m)\bigr)\le \const$. 
So, suppose  that $\sum_{i=1}^n \big(X_i^e(m)-M^e(m)\big)^2>0$.

Using inequality $|\sqrt{1+b}-1|\le |b|$ for $b\ge -1$,
we obtain that
\begin{eqnarray}
\lefteqn{\bigl|g\bigl(\tilde X^e(m+1)\bigr)-g\bigl(\tilde X^e(m)\bigr)\bigr|}\nonumber\\
&=&\bigl[f\bigl(X^e(m)\bigr)\bigr]^{1/2}\biggl|\biggl(1+\frac{f\bigl(X^e(m+1)\bigr)-f\bigl(X^e(m)\bigr)}{f\bigl(X^e(m)\bigr)}\biggr)^{1/2}-1 \biggr|\nonumber\\
&\le&\bigl[f\bigl(X^e(m)\bigr)\bigr]^{1/2}\biggl| \frac{f\bigl(X^e(m+1)\bigr)-f\bigl(X^e(m)\bigr)}{f\bigl(X^e(m)\bigr)}\biggr|\nonumber\\
&=&  \frac{\bigl|f\bigl(X^e(m+1)\bigr)-f\bigl(X^e(m)\bigr)\bigr|}{\bigl[f\bigl(X^e(m)\bigr)\bigr]^{1/2}}\\
&\le& 4, \nonumber
\end{eqnarray}
by Lemma~\ref{lemma_g1}.

Let 
\[
A=\mathcal{M}\cap\{x\in \R^n: \; \max_i |x_i|< a\},
\]
where $a$ is from Lemma~\ref{lemma_g2}.
That is, $A$ is the set of possible configurations of $\tilde X^e$ such that
$\max_{i=1, \ldots, n}|\tilde X^e_i|<a$. Note that the set $A$ is finite.
Let us now prove that
\[
\E\bigl[g\bigl(\tilde X^e(m+1)\bigr)-g\bigl(\tilde X^e(m)\bigr)\mid  \tilde X^e(m)=x \bigr]\le -c_2/\sqrt{2}, 
\]
if $x\notin A$. 
Indeed, if $x\in \mathcal{M}\setminus  A$, as $\sqrt{1+b}\le 1+\frac{b}{2}$ for $b\ge -1$, we get
(using Lemma~\ref{lemma_g2})
\[
\begin{split}
&\E\bigl[g\bigl(\tilde X^e(m+1)\bigr)-g\bigl(\tilde X^e(m)\bigr)\mid  \tilde X^e(m)=x \bigr]\\
&\quad=\sqrt{f(x)}
\, \E\biggl[ \biggl(1+\frac{f\bigl(\tilde X^e(m+1)\bigr)-f\bigl(\tilde X^e(m)\bigr)}{f\bigl(\tilde X^e(m)\bigr)}\biggr)^{1/2}-1
\ \Big| \ \tilde X^e(m)=x \biggr]\\
&\quad\le\frac{\E\bigl[f\bigl(\tilde X^e(m+1)\bigr)-f(x)\mid \tilde X^e(m)=x \bigr]}{2\sqrt{f(x)}}\\
&\quad\le -\frac{c_2}{2}.
\end{split}
\]
Thus, by Theorem~\ref{t:ergodic} the process $\tilde X^e$ is positive recurrent.

For $\tau_A=\inf\{m>0: \tilde X^e(m+k) \in A \} $,
take now
\begin{eqnarray*}
{\mathfrak S}_m=\left\{ \begin{array}{lcl}
g\bigl(\tilde X^e(m)\bigr), &\mbox{ if }& m\le \tau_A,\\
-(m-\tau_A), &\mbox{ if }& m> \tau_A
\end{array}
\right.
\end{eqnarray*}
and apply Theorem~\ref{217} to the sequence
$\{{\mathfrak S}_m\}$.
We have that for any $\delta_1<c_2/2$,  there exist $C$ and $\delta_2$ such that
\[
\P[\tau_A>(1-\delta_1)m \mid \tilde X^e(k)=x'\notin A]<Ce^{-\delta_2 m}.
\]
Note that there exist $k$ and $\delta>0$ such that  for any $y\in A$
\[
\P[ \tilde X^e(m+l)=0\text{ for some }l\le k \mid \tilde X^e(m)=y]\ge \delta.
\]
It is then not difficult to obtain that $\E ( e^{c'\tau}\mid \tilde X^e(k)=x) <\infty$, where 
$\tau=\inf\{m>0: \tilde X^e(m)=0\}$. This proves part (i). 

For the part (ii), fix a routing policy $P$ accordingly to
Definition~\ref{routing_policy} and let  $Z(m)$ be the process obtained using
this routing policy. That is, when $Z(m)=x$, an item that arrives at $S_i$ is
directed to node $s^i_j$ with probability $P^{(i)}_j(x)$, and then
$Z_j(m+1)=Z(m)+1$, $Z_l(m+1)=Z_l(m)$ for $l\ne j$. We will compare this process
to $X^e(m)$ obtained with JSQ routing policy. Let these to processes be defined
at the same probability space and use the same arrivals, but the routing
policies act independently. Analogously to Lemma~\ref{lemmaJSQ}, we get
\begin{align*}
&\E\bigl[f\bigl(X^e(m+1)\bigr)-f\bigl(Z(m+1)\bigr)\mid X^e(m)=Z(m)=x , \; C_i\bigr]\\
&~~~= - 2\sum_{j\ne j_{\min}} P^{(i)}_j(x) \bigl( X^e(i,j,m)-X^e(i,j_{\min},m)\bigr).
\end{align*}
Thus, 
\begin{eqnarray}
\lefteqn{\E\bigl[f\bigl(X^e(m+1)\bigr)-f\bigl(Z(m+1)\bigr)\mid X^e(m)=Z(m)=x \bigr]}\nonumber\\
&=& -2\sum_{i=1}^K  \lambda_i\sum_{j\ne j_{\min}}P^{(i)}_j(x) \bigl(X^e(i,j,m)-X^e(i,j_{\min},m)\bigr)\le 0,
 \end{eqnarray}
which proves part (ii).  Theorem~\ref{time} is proved.
\qed

%

\medskip
\noindent
{\it Proof of Theorem~\ref{erg_alpha}}. 
Let $N_i(t)$ be the number of arrivals at $S_i$ by time $t$. 
Since $N_i(t)$ is a Poisson process with rate $\lambda_i$, a.s.\
$N_i(t) \to \infty$ and $N_i(t)/t\to \lambda_i$ as $t \to \infty$, $i=1, \ldots, n$.

As $\tilde X(t)$ is recurrent,
we have that for almost every realization of the process $\tilde X(t)$ 
there exists an infinite sequence  $t_1, t_2, \ldots$
such that $\tilde X(t_j)=0$ for all $j$.
For these moments $t_j$  we can define
%
\[ 
\alpha_{ik}(t_j) = \frac{\lambda_i N_{ik}(t_j)}{N_i(t_j)},
 \]
where  $N_{ik}(t_j)$ is the number of items arrived at $S_i$ and 
 routed to node $s^i_k$ by time $t_j$.
So,  sending  the proportion
$\frac{\alpha_{ik}(t_j)}{\lambda_i}$ of items arriving at $S_i$
to  $s^i_k$, results in the same number of items at all 
nodes. As the sequence of $\alpha_{ik}(t_j)$ is bounded, we can chose
a subsequence $\alpha_{ik}(t_{j'})\to \alpha_{ik}$, as $t_{j'}\to \infty$.
Evidently, $\alpha_{ik}\ge 0$ and $\sum_{k=1}^{\kappa_i}\alpha_{ik}=\lambda_i$.
Then, as \[
       X_l(t_j)=\frac{1}{n}\sum_{l'=1}^n X_{l'}(t_j)=\frac{1}{n}\sum_{i=1}^K N_i(t_j)    
         \]
 we obtain
\begin{align*}
 \frac{1}{n}&=\frac{X_l(t_j)}{\sum_{i=1}^K N_i(t_j)}\\
& =\frac{1}{ \sum_{i=1}^K N_i(t_j) }\sum_{i=1}^K \sum_{m=1}^{\kappa_i}N_{im}(t_j)\delta_{l, s^i_m}\\
&=\frac{1}{ \sum_{i=1}^K N_i(t_j) }\sum_{i=1}^K N_i(t_j)
\sum_{m=1}^{\kappa_i}\frac{N_{im}(t_j)}{N_i(t_j)}\delta_{l, s^i_m}\\
&=\frac{t_j}{ \sum_{i=1}^K N_i(t_j) }\sum_{i=1}^K \frac{N_i(t_j)}{t_j}
\sum_{m=1}^{\kappa_i}\frac{\alpha_{im}(t_j)}{\lambda_i}\delta_{l, s^i_m}.
\end{align*}
As 
\[
\frac{N_i(t)}{t}\to \lambda_i \text{ and }\frac{t}{ \sum_{i=1}^K N_i(t) }\to \frac{1}{ \sum_{i=1}^K \lambda_i}=1,
\]
and $\alpha_{im}(t_{j'})\to \alpha_{im}$, we see that $\{\alpha_{im}\}$ is indeed a solution of 
\eqref{lin_pr} and
Theorem~\ref{erg_alpha} is proved.
\qed

\subsection{Proof of Theorem \ref{t:noterg}}
\label{s:noterg}
We will use theorems from \cite{R}, therefore let us recall some definitions from there.
A subset $C$ in $\R^n$ is called {\em convex} if 
$(1-\lambda)x+\lambda y \in C$ for every $x \in C$, $y \in C$ and $0<\lambda<1$.
A subset $M$ in $\R^n$ is called an {\em affine set} if 
$(1-\lambda)x+\lambda y \in M$ for every $x \in M$, $y \in M$ and $\lambda \in \R$.
Given any set $A \subset \R^n$ there exists a unique smallest affine set
containing $A$ (namely, the intersection of the collection of the affine sets
$M$ such that $A \subset M$), this set is called {\em affine hull} of $A$ and
is denoted by $\aff A $.
Given a set $A \subset \R^n$ the interior that results when $A$ is
regarded as a subset of its affine hull $\aff A $ is called {\em relative
interior} of $A$ and is denoted by $\ri A$. The closure of $A$ is
denoted by $\cl A$. Note that $\cl(\cl A)=\cl A$ and $\ri (\ri A)=\ri A$;
moreover, if $A$ is convex, then $\cl(\ri A)=\cl A$ (see Theorem~6.3  in~\cite{R}).
If $A$ is convex and $A\ne \varnothing$, then $\ri A \ne \varnothing$ (see
Theorem~6.2 in~\cite{R}).
A set $A$ is called {\em relatively open} if $\ri A=A$.

Let us apply the definitions to our model.
Note that $\Lambda_i\in \R^\kappa_i$ is convex, $i=1,\dotsc, K$.
By $E$ denote the linear transformation that takes a point $p=(p^{(1)},\dotsc,
p^{(K)}) \in \R^{\K}$ to the point $x =(x_1,\dots, x_n) \in \R^n$ such that
\[
x_{\ell}=\sum_{i=1}^K   \sum_{j=1}^{\kappa_i} \lambda_i p^{(i)}_j
\delta_{\ell,s^i_j} \quad \text{for $\ell=1,\dotsc,n$}
\]
where, as before, $\delta_{\ell m}$ is a Kronecker delta. 

\begin{align}
\label{LyD}
\text{ Let } L &:=E(\Lambda_1 \times \dotsb \times \Lambda_K ) \subset \R^n \text{ and } \\ 
D& :=F(L)\subset  {\M}\subset{\R^n}.\nonumber
\end{align}
Since, for $i=1,\dotsc, K$, the set $\Lambda_i$ is convex, we
see that the set $\Lambda_1 \times \dotsb \times \Lambda_K $ is convex
(see Theorem 3.5 in \cite{R}). As
$E$ and $F$ are linear transformations,  the  sets  $L$ and $D$ are
convex (see Theorem 3.4 in \cite{R}).
Since, for $i=1,\dotsc, K$, the set $\Lambda_i$ is compact, 
 the set $\Lambda_1 \times \dotsb \times \Lambda_K $ is also compact.  Since
$E$ and $F$ are linear transformations, and therefore, continuous transformations, we see that the  sets  $L$ and $D$ are
compact. In particular, $D$ is closed, that is, $\cl D=D$.

To translate the condition of Theorem \ref{t:noterg} to the language of convex analysis,
we need the following lemma.
\begin{lm}
\label{lm:riD}
For any parameters of the model $S_1, \ldots, S_K $ and $\lambda_1, \dots, \lambda_n$, the following statements are equivalent:
\begin{enumerate}
\item $0 \in \ri D$;
\item there exists a  positive solution $\alpha_{ij}$ of the system \eqref{lin_pr}.
\end{enumerate}
\end{lm}
{\it Proof}.
Note that  $\ri \Lambda_i$ is the set of points 
$p^{(i)}=(p^{(i)}_1,\ldots, p^{(i)}_{\kappa_i}) \in \R^{\kappa_i}$ such that
\begin{equation}
\label{ri:1}
\begin{cases}
p^{(i)}_j> 0 & \text{for $j=1,\dotsc,\kappa_i$}, \\
\sum_{j=1}^{\kappa_i} p^{(i)}_j=1
\end{cases}
\end{equation}
Moreover, 
\begin{equation}
\label{ri:2}
\ri (\Lambda_1\times \dotsb \times \Lambda_K )= 
(\ri \Lambda_1)\times \dotsb \times (\ri \Lambda_K )
\end{equation}
(see the proof of Corollary 6.6.1 in \cite{R}).
Since $F \circ E$ is a linear transformation, we see that 
\begin{equation*}
\label{ri:3}
F\bigl[E\bigl(\ri (\Lambda_1\times \dotsb \times \Lambda_K )\bigr)\bigr]=
\ri F[E(\Lambda_1\times \dotsb \times \Lambda_K )]
=\ri D
\end{equation*}
(for the first equality see Theorem 6.6. in \cite{R}).


Thus we have proved that $F\bigl[E\bigl((\ri \Lambda_1)\times \dotsb \times (\ri
\Lambda_K )\bigr)\bigr]=\ri D$.  Therefore, $y=F[E(p)]\in \ri D$ if and
only if 
$p \in (\ri \Lambda_1)\times \dotsb \times (\ri \Lambda_K )$.
Recalling \eqref{ri:1} for $\ri \Lambda_i$,
we get that $y=F[E(p)]\in \ri D$ if and only if
\begin{equation}
\label{ri:5}
\begin{cases}
p^{(i)}_j> 0 & \text{for  $j=1,\dotsc,\kappa_i$ and $i=1,\dotsc, K$}, \\
\sum_{j=1}^{\kappa_i} p^{(i)}_j=1 & \text{for $i=1,\dotsc, K$}. 
\end{cases}
\end{equation}

Suppose that item~1 holds, that is, $0 \in \ri D$. Then there exists 
$p\in \R^{\K}$ and $x\in \R^n$ such that $p$ satisfies \eqref{ri:5}, $E(p)=x$ and $F(x)=0$.
Then we have
\[
\sum_{\ell=1}^n x_{\ell} = \sum_{\ell=1}^n 
\sum_{i=1}^K 
 \sum_{j=1}^{\kappa_i} \lambda_i p^{(i)}_j \delta_{l,s^i_j}=
\sum_{i=1}^K 
 \sum_{j=1}^{\kappa_i} \lambda_i p^{(i)}_j=
\sum_{i=1}^K 
  \lambda_i=1.
\]
In the first equation we use the definition of $E$, in the second the fact that 
$\sum_{\ell=1}^n \delta_{l,s^i_j}=1$, in the third the second line from \eqref{ri:5}.
Therefore, using the definition of $F$,  it follows from $F(x)=0$  that $x=E(p)=(\frac{1}{n},\dotsc,\frac{1}{n})$. Then
$\sum_{j=1}^{\kappa_i} \lambda_i p^{(i)}_j \delta_{l,s^i_j}=\frac{1}{n}$ for $\ell=1,\dotsc,n$.
We have proved that if $0 \in \ri D$ then there exists $p\in \R^{\K}$ such that
\begin{equation}
\label{ri:7}
\begin{cases}
p^{(i)}_j> 0 & \text{for $j=1,\dotsc,\kappa_i$ and $i=1,\dotsc, K$}, \\
\sum_{j=1}^{\kappa_i} p^{(i)}_j=1 & \text{for $i=1,\dotsc, K$}, \\
\sum_{j=1}^{\kappa_i} \lambda_i p^{(i)}_j \delta_{l,s^i_j}=\frac{1}{n} &\text{for $\ell=1,\dotsc,n$}.
\end{cases}
\end{equation}
Substituting 
$\frac{\alpha_{ij}}{\lambda_i}$
for $p^{(i)}_j$ in \eqref{ri:7}, we get a positive solution of~\eqref{lin_pr}. Thus item~2 holds.

Now suppose that item~2 holds, that is, there exists $p\in \R^{\K}$ that
satisfies \eqref{ri:7}.  Let us prove that $0\in \ri D$.  Comparing the
first and second line of \eqref{ri:7} with \eqref{ri:5}, we get $F[E(p)] \in
\ri D$. 
Let $x=E(p)$. Then 
\[
x_{\ell} = 
\sum_{j=1}^{\kappa_i} \lambda_i p^{(i)}_j \delta_{l,s^i_j}=
\frac{1}{n}.
\]
In the first equation we use the definition of $E$ and in the second the third line
of \eqref{ri:7}. From the  definition of $F$ it follows that $F(x)=0$. Therefore,
$F[E(p)]=0$. Thus $0=F[E(p)]\in \ri D$ and item~1 holds.  \qed

Let us recall some additional definitions from~\cite{R}.
For $M \subset \R^n$ and $a \in \R^n$, the {\em translate} of $M$ by $a$ is
defined to be set
\[
M+a=\{x+a \mid x \in M\}.
\]
A translate of an affine set is another affine set. An affine set $M$ is {\em parallel} to
an affine set $L$ if $M=L+a$ for some $a$. Each non-empty affine set is
parallel to a unique subspace $L$ (see Theorem 1.2 in \cite{R}). The {\em
dimension} of a non-empty affine set is defined as the dimension of the
subspace parallel to it. An $(n-1)$-dimensional affine set in $\R^n$ is called a
{\em hyperplane}.  
By $\langle \cdot,\cdot\rangle$ denote the inner product in $\R^n$:
$\langle x,y\rangle =\sum_{i=1}^n x_i y_i$.
Given $\beta \in \R$ and a non-zero $b\in \R^n$, the set 
\[
H=\{x \mid  \langle x,b \rangle=\beta\}
\]
is a hyperplane in $\R^n$; moreover, every hyperplane may be represented in
this way, with $b$ and $\beta$ unique up to common non-zero multiple (see
Theorem 1.3 in \cite{R}). 


For any non-zero $b \in \R^n$ and any $\beta \in \R$, the sets
\[
\{ x \mid \langle x, b \rangle  \le \beta\}, \qquad
\{ x \mid \langle x, b \rangle  \ge \beta\}
\]
are called {\em closed half-spaces}.
The sets
\[
\{ x \mid \langle x, b \rangle  < \beta\}, \qquad
\{ x \mid \langle x, b \rangle  > \beta\}
\]
are called {\em open half-spaces}.
The half-spaces depend only on the hyperplane $H=\{x:\langle x, b \rangle =\beta\}$. One may
speak unambiguously, therefore, of the open and closed hyperspaces
corresponding to a given hyperplane.

Let $C_1$ and $C_2$ be non-empty sets in $\R^n$. A hyperplane is said to {\em
separate} $C_1$ and $C_2$ if $C_1$ is contained in one of the closed half
spaces associated with $H$ and $C_2$ lies in the opposite half-space.
It is said that to separate $C_1$ and $C_2$ properly if $C_1$ and $C_2$ are not
{\em both} actually contained in $H$ itself.

Now we are ready to prove Theorem \ref{t:noterg}. By Lemma~\ref{lm:riD}, we have 
$0 \notin \ri D$.

Note that the one point set $\{0\}$ is an affine set, $\ri D$ is a relatively open convex
set and $(\ri D) \cap \{0\}=\varnothing$. Therefore, there exists a hyperplane $H$
containing $0$ such that one of the open half-spaces associated with $H$
contains $\ri D$ (see Theorem 11.2 in \cite{R}).  Since $0\in H$, we see that
$H=\{x:\langle x, b \rangle =0\}$ with some $b\in \R^n$, $b\ne 0$.  Substituting, if it is necessary,
$-b$ for $b$, we see that there is a linear functional $f: \R^n \to
\R$ that sends point $y\in \R^n$ to value $ \langle y, b \rangle $, and if $y \in
\ri D$, then $f(y)>0$. 

Recall that the state space of the Markov chain $\tilde X^e(m)$ is $\mathcal{M}$.
Since $\ri D \subset \M$, we see that there is a point $z \in \mathcal{M}\subset \M$  such that $f(z)>0$. Also, $0 \in \mathcal{M}$ and $f(0)=0$. 
To prove that $\tilde X^e(m)$ is not positive recurrent let us apply 
Theorem~\ref{t:non-ergodic} to the Markov chain $\tilde X^e(m)$
and the function $f$. To apply the theorem, we see that it is enough to check that   
\begin{equation}
\label{eq:tildex}
\E\bigl[f\bigl(\tilde X^e(m+1)\bigr)-f\bigl(\tilde X^e(m)\bigr)\mid \tilde X^e(m)=z\bigr] \ge 0
\quad \text{for any $z\in \mathcal{M}$}
\end{equation}
To prove \eqref{eq:tildex} it is enough to prove 
\begin{equation}
\label{eq:x}
\E\bigl[f\bigl(F[X^e(m+1)]\bigr)-f\bigl(F[X^e(m)]\bigr)\mid X^e(m)=x\bigr] \ge 0 \quad \forall x\in \N^n
\end{equation}

To prove \eqref{eq:x}, we need some notation.
For $i=1,\dots, K$, denote by $e^{(i)}_j$ the $j$-th coordinate vector in
$\R^{\kappa_i}$.
By $T$ denote the linear transformation that takes a point $p=(p^{(1)},\dotsc,
p^{(K)}) \in \R^{\K}$ to the point $x =(x_1,\dots, x_n) \in \R^n$ such that
\[
x_{\ell}=\sum_{i=1}^K   \sum_{j=1}^{\kappa_i} p^{(i)}_j
\delta_{\ell,s^i_j} \quad \text{for $\ell=1,\dotsc,n$}
\]
In particular, $T$ takes the point $e^{(i)}_j$ to the point $x$ such that
$x_{\ell}=\delta_{\ell,s^i_j}$, for $l=1,\dots,n$.

Let us prove \eqref{eq:x}.
Take any $x \in \N^n$.
Recall that, for routing policy $P$,
we have $p=P(x)\in \Lambda_1\times \dotsb \times \Lambda_K $.
Moreover, 
\begin{align*}
&\E\bigl[f\bigl(F[X^e(m+1)]\bigr)-f\bigl(F[X^e(m)]\bigr)\mid X^e(m)=x\bigr]\\
&\quad 
=
\sum_{i=1}^K  \sum_{j=1}^{\kappa_i}
\lambda_i p^{(i)}_j\bigl\{f\bigl(F[x+T(e^{(i)}_j)]\bigr)-f\bigl(F[x]\bigr)\bigr\}\\
&\quad 
=
\sum_{i=1}^K  \sum_{j=1}^{\kappa_i}
\lambda_i p^{(i)}_jf\bigl(F[T(e^{(i)}_j)]\bigr)\\
&\quad 
=
f\biggl(F\Bigl[\sum_{i=1}^K  \sum_{j=1}^{\kappa_i}
\lambda_i p^{(i)}_jT(e^{(i)}_j)\Bigr]\biggr)
=
f\bigl(F[E(p)]\bigr) \ge 0.
\end{align*}
In the second and third equalities we use that $f\circ F$ is 
a linear functional. Let us check the last inequality. 
We have 
\[
F[E(p)] \in F[E(\Lambda_1\times \dotsb \times \Lambda_K )]=D.
\]
Since $D$ is closed and convex, we see that $\cl (\ri
D)=\cl D=D$ (see the properties of operations $\ri$ and $\cl$ in the beginning of
Section~\ref{s:noterg}).
Note that $f(y)>0$ for $y \in \ri D$ and linear functional $f$ is continuous,
therefore, $f(y)\ge 0$ for $y \in \cl(\ri D)=D$.  

Thus all conditions of Theorem~\ref{t:non-ergodic} are satisfied and Theorem~\ref{t:noterg} is proved.
\qed


\section*{Acknowledgements}
The work of M.M. was partly supported by CNPq (grant 450787/2008--7), FAPESP (grants 2007/50459--9, 2004/03056--8).

The work of V.S. was  supported by FAPERJ (grants E-26/170.008/2008
and E-26/110.982/2008) and CNPq (grants 471891/2006-1, 309397/2008-1
and 471946/2008-7).

The work of M.V. was partly supported by CNPq (grants 304561/2006--1, 301455/2009--0,  472431/2009--9, 471925/2006--3) and FAPESP (thematic grants 09/52379--8, 2004/07276--2).

We thank Iain MacPhee who has read  the manuscript very carefully and helped us to improve the presentation.
Also, the authors thank the anonymous referee for valuable comments
 and suggestions.

\end{document}